\documentclass[12pt]{amsart}
\usepackage{amsmath,mathrsfs}
\usepackage{amssymb}
\usepackage{latexsym}
\usepackage{amsfonts}
\newtheorem{Pa}{Paper}[section]
\newtheorem{Tm}[Pa]{{\bf Theorem}}
\newtheorem{La}[Pa]{{\bf Lemma}}
\newtheorem{Cy}[Pa]{{\bf Corollary}}
\newtheorem{Rk}[Pa]{{\bf Remark}}
\newtheorem{Pn}[Pa]{{\bf Proposition}}

\newtheorem{Dn}[Pa]{{\bf Definition}}

\def\N{\mathbb N}
\def\C{\mathbb C}
\def\z{\zeta}
\begin{document}
\title[Linear state space theory and the white noise space]
{Linear State space theory in the white noise space setting}
\author[D. Alpay]{Daniel Alpay}
\address{(DA) Department of mathematics,
Ben-Gurion University of the Negev, P.O. Box
653, Beer-Sheva 84105, Israel}
\email{dany@math.bgu.ac.il}
\author[D. Levanony]{David Levanony}
\address{(DL) Department of electrical
engineering, Ben-Gurion University of the Negev,
P.O. Box 653, Beer-Sheva 84105, Israel}
\email{levanony@ee.bgu.ac.il}
\author[A. Pinhas]{Ariel Pinhas}
\address{(AP) Department of mathematics,
Ben-Gurion University of the Negev, P.O. Box 653, Beer-Sheva
84105, Israel}
\email{arielp@bgu.ac.il}
\thanks{D. Alpay thanks the
Earl Katz family for endowing the chair which supported his
research. This research is part of the European Science
Foundation Networking Program HCAA, and was supported in part by
the Israel Science Foundation grant 1023/07}
%\begin{document}
\date{}

\begin{abstract}
We study state space equations within the white noise space
setting. A commutative ring of power series in a countable number
of variables plays an important role. Transfer functions are
rational functions with coefficients in this commutative ring,
and are characterized in a number of ways. A major feature in our
approach is the observation that key characteristics of a linear,
time invariant, stochastic system are determined by the
corresponding characteristics associated with the deterministic
part of the system,
namely its average behavior.\\

\end{abstract}

%\maketitle
\maketitle
\tableofcontents

\section{Introduction}
\setcounter{equation}{0}
In a preceding paper, see \cite{al_acap}, the first two authors
began a study of linear stochastic systems within the framework
of the white noise space. There, the emphasis was on stability
theorems associated with convolution systems. The present paper
is concerned with state space theory. Specifically, we study
systems defined by state space equations when randomness is
allowed in the matrices defining these equations.

\begin{Rk}
While the discussion to follow is restricted to discrete time,
most results apply to continuous time in an obvious way.
Continuous time is not explicitly pursued in this work.
\end{Rk}

To set the problem in perspective and provide
motivation, we first recall some well known facts from linear
system theory. There, state space equations of the form
\begin{equation}
\label{eq:stp}
\begin{split} x_{n+1}&=Ax_n+Bu_n,\quad n=0,1,\ldots\\
                      y_n&=Cx_n+Du_n,
\end{split}
\end{equation}
play an important role. In \eqref{eq:stp},  $A\in{\C}^{N\times
N}$, $B\in{\C}^{N\times q}$, $C\in{\C}^{p\times N}$,
$D\in{\C}^{p\times q}$, the states $x_n$ take values in $\C^N$,
the inputs $u_n$ in $\C^q$ and the outputs $y_n$ in $\C^p$. Taking the Z
transform, assuming $x_0=0$, \eqref{eq:stp} leads to
\begin{equation}
\begin{split}
\mathscr X(\z)&=\z A\mathscr X(\z)+\z B\mathscr U(\z)\\
\mathscr Y(\z)&=C\mathscr X(\z)+D\mathscr U(\z),
\end{split}
\end{equation}
where the Z transform variable is denoted by $\z$, so that
\begin{equation}
\mathscr Y(\zeta)=\mathscr H(\zeta)\mathscr U(\zeta),
\end{equation}
where
\[
\mathscr Y(\zeta)=\sum_{n=0}^\infty y_n\zeta^n,\quad \mathscr
U(\zeta)=\sum_{n=0}^\infty u_n\zeta^n,\quad \mathscr
X(\zeta)=\sum_{n=0}^\infty x_n\zeta^n,
\]
and
\begin{equation}
\label{eq:tf1}
\mathscr H(\zeta)=D+\zeta C(I_N-\zeta A)^{-1}B.
\end{equation}
The function $\mathscr H(\z)$ is called the transfer function of
the system defined by \eqref{eq:stp}. It is also possible to
replace ${\mathbb C}$ by a field ${\mathfrak K}$ over the complex
numbers. Matrices $A,B,C$ and $D$ then have their entries in
${\mathfrak K}$, and $\mathscr H(\z)$ makes sense for all $\zeta$
such that $(I_N-\z A)$ is invertible. The case where the
components of $A,B,C$ and $D$ belong to some commutative ring is
of special interest. See for instance \cite{MR525830},
\cite{MR0427298}, \cite{SSR}, \cite{MR1648312}, \cite{MR0286517},
\cite{MR839186}. This ring is often assumed Noetherian, to enable
to formulate results. Even when the ring contains the complex
numbers, formula \eqref{eq:tf1} does not make sense in general
because elements in the ring need not be invertible in the ring.
It will make sense in a normed ring for $\z$ small enough, as is
seen using the power expansion
\begin{equation}
\label{eq:nr} (I_N-\z A)^{-1}=\sum_{n=0}^\infty \z^n A^n.
\end{equation}
The ring $\mathfrak R$ defined below is not a normed ring, but
still it will be possible to define expansions of the form
\eqref{eq:nr} in it.\\

As we have explained in our previous paper \cite{al_acap}, a
Gaussian input into a linear system with nonrandom coefficients,
will result in a Gaussian output. In the present series of
papers, and in particular in the present work, we aim to model
linear Gaussian input-output relations when the underlying linear
system is random. Here we allow Gaussian inputs and randomness in
the matrices $A,B,C$ and $D$, in such a way that the outputs
remain Gaussian. While indeed a Gaussian input into a linear
system with random coefficients cannot be expected to result in a
Gaussian output, we will use the white noise space setting (see
\cite{MR1244577}, \cite{MR1408433} and Section \ref{2} below), and
replace the pointwise product with the Wick product, enabling
Gaussian input-output relations when the underlying system has
random coefficients. This framework will preserve Gaussian
input-output relation, while allowing uncertainty in the form of
randomness in the linear system under study. Such a setting may
prove useful so as to study a linear state space system with {\it
nonrandom} uncertainties, a system that indeed maintains Gaussian
input-output relation. This, by utilizing the Bayesian embedding
approach to solve problems associated with a system subjected to a
deterministic uncertainty, through solutions to corresponding
problems associated  with systems  with random
uncertainties, see e.g. \cite{MR1931659}.\\

In the white noise space setting, the space of complex numbers
$\mathbb C$ is replaced by a space of stochastic distributions
called the Kondratiev space, denoted by $S_{-1}$. This space
contains ${\mathbb C}$, is the inductive limit of a certain
family of Hilbert spaces (see \eqref{michelle} below), and is
nuclear; see \cite[Definition 2.3.2 (b) p. 30, Lemma 2.8.2 p.
74]{MR1408433}. A key element in our formulation is a product
defined on ${S}_{-1}$, namely the Wick product, denoted by
$u\lozenge v$, which reduces to multiplication by a constant when
one of the elements $u$ or $v$ is non random. We thus replace the
equations \eqref{eq:stp} by
\begin{equation}
\label{eq:sp2}
\begin{split}
x_{n+1}&=A\lozenge x_n+B\lozenge u_n\\
y_n&=C\lozenge x_n+D\lozenge u_n
\end{split}
\end{equation}
where $A\in({S_{-1}})^{N\times N}$,  $B\in({S_{-1}})^{N\times q}$,
$C\in({S_{-1}})^{p\times N}$, and  $D\in({S_{-1}})^{p\times q}$.
The states $x_n$ take values in $(S_{-1})^N$, the input $u_n$ in
$(S_{-1})^q$ the output in $(S_{-1})^p$.\\

A fundamental tool in white noise analysis is the Hermite
transform
\[
F\mapsto{\mathbf I}(F)
\]
(see below), which associates to every
element in $S_{-1}$, a power series in a countable number of
complex variables,
\begin{equation}
z=(z_1,z_2,z_3,\ldots)
\end{equation}
and transforms the Wick product into a point-wise product:
\begin{equation}
\label{hermite product} \textbf{I}(F\diamond G)(z)=
(\textbf{I}(F)(z))(\textbf{I}(G)(z)),\quad \forall F,G\in S_{-1}.
\end{equation}
The image of the Kontradiev space under the Hermite transform is a
commutative ring without divisors of zeros (that is, a domain),
which we will denote by ${\mathfrak R}$. It is not Noetherian, so
most results in system theory on commutative rings cannot be
applied. Still, it has a very important property, which allows us
to proceed. An $\mathbf F\in{\mathfrak R}$ is invertible in
${\mathfrak R}$ if and only if its constant coefficient is non
zero (recall that $\mathbf F$ is a power series). More generally,
for $p\in {\mathbb N}$, an $\mathbf F\in{\mathfrak R}^{p\times
p}$ will be invertible in ${\mathfrak R}^{p\times p}$ if and only
if the matrix $\mathbf F(0)$ (which belongs to ${\mathbb
C}^{p\times p}$) is invertible; see Theorem \ref{tm:inv} below.
This theorem follows from a non trivial result on the
characterization of the range of the Hermite transform, given in
\cite[Theorem 2.6.11,
p. 62]{MR1408433}.\\

We now take the Z transform  and the Hermite transform of
\eqref{eq:sp2}. The Z transform of the series $({\mathbf
I}(u_n))_{n=0,1,\ldots}$ is denoted by $\mathscr U(\z,z)$
\[
\mathscr U(\z,z)=\sum_{n=0}^\infty \zeta^n({\mathbf I}(u_n))(z),
\]
and similarly for $\mathscr Y(\z,z)$ and $\mathscr X(\z,z)$. We
obtain:
\begin{equation}
\label{1.9}
\begin{split}
(I_N-\zeta \mathbf A(z))\mathscr X(\zeta,z)&=\mathbf B(z)
\mathscr U(\zeta,z),\\
\mathscr Y(\zeta, z)&=\mathbf C(z)\mathscr X(\zeta, z)+\mathbf
D(z)\mathscr U(\zeta, z),
\end{split}
\end{equation}
where
\[
\mathbf A(z)={\mathbf I}(A)(z),\,\,\mathbf B(z)={\mathbf
I}(B)(z),\,\,\mathbf C(z)={\mathbf I}(C)(z),\,\,\mathbf
D(z)={\mathbf I}(D)(z).
\]
Since $\mathbf A(z)$ is bounded in a certain neighborhood of the
origin (see Theorem \ref{tm:inv} below), the matrix $(I_N-\zeta
\mathbf A(z))$ is invertible for $\zeta$ small enough, and we can
write
\[
\mathscr Y(\z,z)={\mathscr H}(\z,z)\mathscr U(\z,z),
\]
where
\begin{equation}
\label{real}
{\mathscr H}(\zeta, z)=\mathbf D(z)+\zeta \mathbf
C(z)(I_N-\zeta \mathbf A(z))^{-1}\mathbf B(z),
\end{equation}
is the transfer function of the system defined by the
equations \eqref{eq:sp2}.\\

\begin{Rk}
When we set $z=0$ in \eqref{1.9} we retrieve \eqref{eq:tf1}, that
is, we are back within the deterministic setting.
\end{Rk}

We now view \eqref{real} as a ${\mathfrak R}$-valued function.
Consider $\z$ such that
\begin{equation}
\label{eq:spec} \det(I_N-\z \mathbf A(0))\not =0.
\end{equation}
It follows from Theorem \ref{tm:facto} below that $(I_N-\zeta
\mathbf A)$ is invertible in ${\mathfrak R}$ for such $\z$.
Therefore the ${\mathfrak R}$-valued function $\mathscr H(\z)$
given by
\[
(\mathscr H(\z))(z)={\mathscr H}(\zeta, z),
\]
that is,
\begin{equation}
\label{eq:real2}
\mathscr H(\z)= \mathbf D+\zeta \mathbf
C(I_N-\zeta \mathbf A)^{-1}\mathbf B
\end{equation}
is well defined for $\zeta$ satisfying \eqref{eq:spec}.\\

Functions of the form \eqref{eq:real2} will be called {\it
rational functions associated with the white noise space}. We
note that in \cite{alpay-2008} another approach to rational
functions, with emphasis on rationality with respect to a finite
number of the variables $z_k$ is considered. The purpose of
this paper is to give a number of equivalent characterizations of
functions of the form \eqref{eq:real2} and to study
the notions of controllability, observability and minimality
in the setting of the ring $\mathfrak R$.\\

A major non-trivial feature of this work is the observation that
key characteristics of a linear, time invariant, stochastic system
(e.g. invertibility), are determined by the corresponding
characteristics associated with the deterministic part of the
system under study, namely its average behavior. For instance, a
realization of the perturbed system \eqref{eq:real2} will be
observable (see Definition \ref{dn:obs} below) if the
corresponding realization of the unperturbed system, namely with
$z=0$ is observable. See Theorem \ref{tm:min}.\\

The paper consists of six sections besides the introduction, and
its outline is as follows. We review, as already mentioned, white
noise space theory in Section \ref{2}. We study the ring
${\mathfrak R}$ in Section \ref{3}. Equivalent characterizations
of rational functions are given in Section \ref{4}. Observable
pairs are studied in Section \ref{5}. In Section 6 we consider
controllable pairs, and briefly discuss minimal realizations.
%Section 7 briefly discusses a subring of $\mathfrak R$.
The last
section considers the case where the functions take values in one
of the Hilbert spaces
which make ${\mathfrak R}$. \\

\section{A brief survey of white noise space analysis}
\label{2} \setcounter{equation}{0}
The starting point to construct the white noise space is the
Schwartz space $\mathcal S$ of {\sl real-valued} smooth functions
which, together with their derivatives, decrease rapidly to zero
at infinity. For $s\in{\mathcal S}$, let $\|s\|$ denote its
${\mathbf L}_2({\mathbb R})$ norm. The function
\[
K(s_1-s_2)=e^{-\frac{\|s_1-s_2\|^2}{2}}
\]
is positive (in the sense of reproducing kernels) for $s_1,s_2$
running in ${\mathcal S}$. The space is nuclear. By an extension
of Bochner's theorem to nuclear spaces due to Minlos (see
\cite{MR0154317}, \cite[Th\'eor\`eme 3, p. 311]{MR35:7123}),
there exists a probability measure $P$ on ${\mathcal S}^\prime$
such that
\[
K(s)=\int_{{\mathcal S}^\prime}e^{-i\langle
s^\prime,s\rangle}dP(s^\prime),
\]
where we have denoted by $\langle s^\prime,s\rangle$ the duality
between ${\mathcal S}$ and ${\mathcal S}^\prime$. The real Hilbert
space ${\mathbf L}_2({\mathcal S}^\prime, {\mathcal F},dP)$, where
${\mathcal F}$ is the Borelian $\sigma$-algebra, is called the
white noise space. We will denote it by $\mathcal W$, and its
elements by $\omega$, by setting $\Omega={\mathcal S}^\prime$.\\

Among all orthogonal Hilbert bases of the white noise space, one
plays a special role. It is constructed in terms of Hermite
functions and its elements are denoted by $H_\alpha$, where the
index $\alpha$ runs through the set $\ell$ of sequences
$(\alpha_1,\alpha_2,\ldots)$, whose entries are in
\[
\mathbb N_0=\left\{
0,1,2,3,\ldots\right\},
\]
and $\alpha_k\not =0$ for only but a finite number of indices $k$.
Furthermore, with the multi-index notation
\[
\alpha!=\alpha_1!\alpha_2!\cdots,
\]
we have
\begin{equation}
\label{eq:fock1}
\|H_\alpha\|_{\mathcal W}^2=\alpha!.
\end{equation}
In view of \eqref{eq:fock1}, the map
\[
H_\alpha\mapsto
z^\alpha=z_1^{\alpha_1}z_2^{\alpha_2}z_3^{\alpha_2}\cdots
\]
extends to a unitary map between ${\mathcal W}$ and the
reproducing kernel Hilbert space with reproducing kernel
\[
k(z,w)=e^{\langle
z,w\rangle_{\ell_2}}=\sum_{\alpha\in\ell}\frac{z^\alpha
w^{*\alpha}}{\alpha!},
\]
where $z,w$ run through $\ell_2$.\\

The Wick product in ${\mathcal W}$ is defined by the formula
\[
H\alpha\lozenge H_\beta=H_{\alpha+\beta},\quad \alpha,
\beta\in\ell,
\]
and the Hermite transform is defined through linearity as
\[
({\mathbf I}(H_\alpha))(z)=z^\alpha.
\]
The space ${\mathcal W}$ is too small to be stable under the Wick
product, and one defines the Kondratiev space $S_{-1}$, within
which the Wick product is stable. More precisely, $S_{-1}$ is
%the dual of
a nuclear space, and is defined as the inductive limit of the
increasing family of Hilbert spaces  ${\mathcal H}_{k},
k=1,2,\ldots$ of formal series $\sum_{\alpha\in\ell}f_\alpha
H_\alpha$ such that
\begin{equation}
\label{michelle}
\|f\|_{k}\stackrel{\rm def.}{=}
 \left(\sum_{\alpha\in\ell}|f_\alpha|^2
(2{\mathbb N})^{-k\alpha}\right)^{1/2}<\infty,
\end{equation}
where, for $\beta\in\ell$,
\[
(2{\mathbb N})^\beta=2^{\beta_1}(2\times 2)^{\beta_2} (2\times
3)^{\beta_3}\cdots.
\]
That the Wick product is stable within $S_{-1}$ is made more
precise by V\r{a}ge's inequality (see \cite[Proposition 3.3.2,
  p. 118]{MR1408433}), which we now recall.
Let $l$ and $k$ be natural numbers such that $k>l+1$.
Let $h\in {\mathcal H}_{l}$ and $u\in
{\mathcal H}_{k}$. Then,
\begin{equation}
\label{vage} \|h\lozenge u\|_{k}\le
A(k-l)\|h\|_{l}\|u\|_{k},
\end{equation}
where
\begin{equation}
\label{vage111}
A(k-l)=\sum_{\alpha\in\ell}(2{\mathbb
N})^{(l-k)\alpha}.
\end{equation}
For a proof that $A(k-l)$ is finite, see
\cite[Proposition 2.3.3, p. 31]{MR1408433}.\\

The series $\sum_{\alpha\in\ell} s_\alpha z^\alpha$ will be said
to be convergent at $z$ if
\[
\sum_{\alpha\in\ell}|s_\alpha||z|^\alpha<\infty,
\]
that is, if it is absolutely convergent, see \cite[p.
60]{MR1408433}. The following easy lemma will be used below.
\begin{La}
Assume that $f(z)=\sum_{\alpha\in\ell}f_\alpha z^\alpha$ and
$g(z)=\sum_{\alpha\in\ell}g_\alpha z^\alpha$ are absolutely
convergent power series at $z$. Then
\begin{equation}
\label{eq:ineq3} |f(z)g(z)|\le \sum_{\gamma\in\ell}|z|^\gamma
\cdot\big| \sum_{\substack{\alpha+\beta=\gamma\\
\alpha,\beta\in\ell}}f_\alpha g_\beta\big|\le
(\sum_{\alpha\in\ell}|f_\alpha|
|z|^\alpha)(\sum_{\alpha\in\ell}|g_\alpha| |z|^\alpha),
\end{equation}
and in particular, the product $fg$ is an absolutely convergent
power series at $z$, and it holds that
\[
f(z)g(z)=\sum_{\gamma\in\ell}z^\gamma\left(\sum_{\substack{\alpha+
\beta=\gamma\\
\alpha,\beta\in\ell}} f_\alpha g_\beta\right).
\]
Furthermore, $(f(z))^n$ is an absolutely convergent power series
at $z$ for all $n\in{\mathbb N}$.
\end{La}
{\bf Proof:} The power series
\[
\sum_{\gamma\in\ell}z^\gamma\left( \sum_{\substack{\alpha+\beta=\gamma\\
\alpha,\beta\in\ell}}f_\alpha g_\beta\right)
\]
is absolutely convergent since
\[
\begin{split}
\sum_{\gamma\in\ell}|z|^\gamma\cdot\big|
\sum_{\substack{\alpha+\beta=\gamma\\
\alpha,\beta\in\ell}}f_\alpha g_\beta\big|&\le
\sum_{\gamma\in\ell}|z|^\gamma\left(
\sum_{\substack{\alpha+\beta=\gamma\\
\alpha,\beta\in\ell}}|f_\alpha| \cdot|g_\beta|\right)\\
&=(\sum_{\alpha\in\ell}|f_\alpha|
|z|^\alpha)(\sum_{\alpha\in\ell}|g_\alpha| |z|^\alpha).
\end{split}
\]
\mbox{}\qed\mbox{}\\

\section{The ring ${\mathbf I}(S_{-1})$}
\setcounter{equation}{0} \label{3} Consider the image ${\mathfrak
R}\stackrel{\rm def.}{=}{\mathbf I}(S_{-1})$ under the Hermite
transform of the Kondratiev space. This is a space of power
series which has been characterized in \cite[Theorem 2.6.11, p.
62]{MR1408433}. In that statement, $(\C^\N)_c$ denotes the space
of finite sequences of complex numbers indexed by the integers,
and the set $K_q(\delta)$ is defined by
\begin{equation}
\label{eq:k_delta} K_q(\delta)= \{z\in\C^\N:
\sum_{\substack{\alpha\in\ell\\
\alpha\not=(0,0,\ldots)}} \left|z^\alpha\right|^2
(2\N)^{q\alpha}<\delta^2\}.
\end{equation}
Note that $\alpha=(0,0,\ldots)$ is excluded from the sum. See
\cite[Definition 2.6.4, p. 59]{MR1408433}.

\begin{Tm}\cite[Theorem 2.6.11, p. 62]{MR1408433}
\begin{enumerate}
\item If $F(\omega)=\sum_\alpha
 a_\alpha H_\alpha (\omega)\in S_{-1}$,
then there exist $q<\infty,~M_q<\infty$ such that
\begin{equation}
\label{eq:ineq}
    \left|\mathbf{I}(F)(z)\right|\leq
    \sum_{\alpha\in\ell} \left|a_\alpha\right|\left|
    z^\alpha\right|\leq M_q\left(\sum_{\alpha\in\ell}
    (2\N)^{q\alpha}\left|z^\alpha\right|^2
    \right)^\frac{1}{2}
\end{equation}
for all $z\in (\C^\N)_c$.
    In particular, $\mathbf I(F)$ is a
    bounded analytic function on $K_q(\delta)$ for
    all $\delta<\infty$.\\
    \item Conversely, suppose
    $g(z)=\sum_\alpha b_\alpha z^\alpha$ is a
    given power series of $z\in (\C^\N)_c$ with
     $b_\alpha\in\C$, with
     $q<\infty$ and $\delta >0$ such that
      $g(z)$ is absolutely convergent for
       $z\in K_q(\delta)$ and
    \[
    \sup_{z\in K_q(\delta)} \left|g(z)\right|<\infty.
    \]
    Then there exists a unique $G\in S_{-1}$
    such that $\mathbf I(G)=g$, namely
    \[
    G(\omega)=\sum_{\alpha\in\ell} b_\alpha H_\alpha (\omega).
    \]
\end{enumerate}
\label{tm:inv}
\end{Tm}

A characterization of convergent sequences in $S_{-1}$ is given
in the following theorem proved in \cite{MR1408433}. It will be
used in particular in the proof of Proposition \ref{pn:obs}.

\begin{Tm}
\label{300408}
\cite[Theorem 2.8.1, p. 74]{MR1408433} A sequence
of elements $F^{(n)}$ in the Kondratiev space $S_{-1}$ converges
to $F\in S_{-1}$ if there exist $\delta>0$ and $q<\infty$ such
that ${\bf I}(F^{(n)})$ converges to ${\bf I}(F)$ pointwise
boundedly, or equivalently, uniformly, in $K_q(\delta)$.
\end{Tm}

The main result of this section is:
\begin{Tm}
$\mathfrak R$ is a commutative ring, which contains ${\mathbb C}$
and has no divisors of zero. Furthermore, let
$x(t)=\sum_{n=0}^\infty x_nt^n$ be a power series, with strictly
positive radius of convergence. Let $p\in{\mathbb N}$. Then for
every $\mathbf r\in\mathfrak R^{p\times p}$ such that $\mathbf
r(0)=0_{p\times p}$, the series
\[
(x(\mathbf r))(z)=\sum_{n=0}^\infty x_n(\mathbf r(z))^n
\]
converges to a limit in $\mathfrak R^{p\times p}$. If
$y(t)=\sum_{n=0}^\infty y_nt^n$ is another such power series, then
\begin{equation}
\label{xy}
(xy)(\mathbf r)=x(\mathbf r)y(\mathbf r),
\quad\forall \mathbf r\in\mathfrak R.
\end{equation}
In particular, an element $\mathbf s$ is invertible in $\mathfrak
R^{p\times p}$ if and only if $\mathbf s(0)$ is invertible.
\label{tm:facto}
\end{Tm}
{\bf Proof:} To simplify the notation we give a proof for $p=1$.
The fact that we have a ring follows from the formula
\eqref{hermite product}. The way to prove the second claim is
 to use Theorem \ref{tm:inv} to show that the {\it
a-priori} formal power series
\[
\sum_{n=0}^\infty x_n(\mathbf r(z))^n
\]
is in fact the image under the Hermite transform of an element in
$S_{-1}$. Since $\mathbf r$ is the image of an element of $S_{-1}$
under the Hermite transform, it satisfies \eqref{eq:ineq} for some
$q\in{\mathbb N}$ and a constant $M_q>0$. Let $r_x$ be the radius
of convergence of the power series defining $x$ (and similarly for
$r_y$ below). We choose $\delta$ such that
\[
M_q\delta\stackrel{\rm def.}{=}\rho<r_x,
\]
Then, by \eqref{eq:ineq}, we have for $z\in K_q(\delta)$,
\[
|\mathbf r(z)|\le\rho,
\]
hence
\[
|\sum_{n=1}^\infty x_n(\mathbf r(z))^n|\le \sum_{n=1}^\infty
|x_n|\rho^n,\quad z\in K_q(\delta).
\]
We conclude the proof by using Theorem \ref{tm:inv}. We first
prove \eqref{xy}. By the preceding arguments we know that
$x(\mathbf r)$, $y(\mathbf r)$, and $(xy)(\mathbf r)$,
are well defined. On the other hand, for $|\mathbf
r(z)|<\min{(r_x,r_y)}$
we have:
\[
\begin{split}
(x(\mathbf r(z))(y(\mathbf r(z))&=\sum_{n=0}^\infty(\sum_{p=0}^n
  x_py_{n-p}) (\mathbf r(z))^n\\
&=((xy)(\mathbf r))(z).
\end{split}
\]

We now turn to the last statement. Assume that $\mathbf s$ is
invertible in $\mathfrak R$, and let $\mathbf u\in\mathbf R$ be
such that $\mathbf s\mathbf u=1$. Then, in particular, $\mathbf
s(0)\mathbf u(0)=1$, so $\mathbf s(0)\not=0$. Conversely, we can
assume without loss of generality, that $\mathbf s(0)=1$. It
suffices then to take in \eqref{xy} $x(t)=1-t$, $y(t)=(1-t)^{-1}$,
and $\mathbf r=\mathbf u-1$. (Note that $\mathbf r(0)=0$.)
\mbox{}\qed\mbox{}\\

\section{Rational functions}
\label{4} Let $f(\zeta, z)=\sum_{n=0}^\infty
f_n(z)\zeta^n\in{\mathfrak R}^{p\times p}(\zeta)$ be a power
series with coefficients in ${\mathfrak R}^{p\times q}$. Define
\[
R_0f(\zeta,z)=\dfrac{f(\zeta,z)-f(0,z)}{\z}.
\]
\begin{Tm}
Let $\mathcal H(\zeta)=\sum_{n=0}^\infty \mathbf
f_n\zeta^n\in{\mathfrak R}^{p\times q}((\zeta)))$ be a formal
power series. Then the
following are equivalent:\\
$(1)$ Components of $\mathcal H$ are obtained by adding,
multiplying and dividing polynomials of ${\mathfrak R}[\zeta]$,
with division being performed only when the constant coefficient
is invertible in ${\mathfrak
R}$.\\
$(2)$ $\mathcal H$ admits a realization in the form of
\eqref{eq:real2}, with
coefficients matrices having entries in ${\mathfrak R}$.\\
$(3)$  The formal power series converges in a neighborhood of the
origin, and there exists a finite number $M$ such that for every
$n\ge M$ the function $R_0^n\mathcal H$ is a linear combination of
$R_0\mathcal H,\ldots R_0^{M-1}\mathcal H$ with coefficients in
$\mathfrak R$.
\end{Tm}
{\bf Proof:} We first note that elements of the form
\eqref{eq:real2} are convergent power series in ${\mathbb R}$ and
not only formal power series, as follows from Theorem
\ref{tm:facto}. Elements $\mathcal H \in{\mathfrak R}^{p\times
q}(\zeta)$  of the form
\begin{equation}
\label{eq:pol}
\mathcal H(\zeta)=\mathbf D\quad{\rm or}\quad \mathcal H(\zeta)=\zeta
\mathbf C,
\end{equation}
where $\mathbf C, \mathbf D\in{\mathfrak R}^{p\times q}$ are clearly
in the form
\eqref{eq:real2}. Furthermore, as is well known, if $p=q$ and
\[
\mathscr H(\z)= \mathbf D+\zeta \mathbf
C(I_N-\zeta \mathbf A)^{-1}\mathbf B\in{\mathfrak R}^{p\times p}(\zeta)
\]
with $\mathbf D$ invertible, then we have:
\[
(\mathscr H(\z))^{-1}=
\mathbf D^{-1}-\zeta \mathbf D^{-1}\mathbf
C(I_N-\zeta \mathbf A^\times)^{-1}\mathbf B \mathbf D^{-1},
\]
where
\[
\mathbf A^\times=\mathbf A-\mathbf B\mathbf D^{-1}\mathbf C.
\]
Furthermore,
if
\[
\mathscr H_1(\z)= \mathbf D_1+\zeta \mathbf
C_1(I_{N_1}-
\zeta \mathbf A_1)^{-1}\mathbf B_1\in{\mathfrak R}^{p_1\times s}(\zeta)
\]
and
\[
\mathscr H_2(\z)= \mathbf D_2+\zeta \mathbf
C_2(I_{N_2}-
\zeta \mathbf A_2)^{-1}\mathbf B_2\in{\mathfrak R}^{s\times
  q_2}(\zeta),
\]
then
\[
(\mathscr H_1\mathscr H_2)(\z)=\mathbf D+\zeta \mathbf
C(I_N-\zeta \mathbf A)^{-1}\mathbf B,\quad N=N_1+N_2,
\]
with $\mathbf D=\mathbf D_1\mathbf D_2$ and
\[
\mathbf A=\begin{pmatrix} \mathbf A_1&\mathbf B_1\mathbf C_2\\
0&\mathbf A_2\end{pmatrix},\quad \mathbf B=\begin{pmatrix}\mathbf
B_1\mathbf D_2\\ \mathbf B_2\end{pmatrix},\quad\mathbf
C=\begin{pmatrix}
\mathbf C_1&\mathbf D_1\mathbf C_2\end{pmatrix}.
\]
A sum of matrices is a special case of a product, as follows from
the formula
\[
\mathbf M_1+\mathbf M_2=\begin{pmatrix} \mathbf
M_1&I_p\end{pmatrix}
\begin{pmatrix} I_q\\ \mathbf M_2\end{pmatrix},
\]
where $\mathbf M_1$ and $\mathbf M_2\in{\mathfrak R}^{p\times
q}$.\\

See for instance \cite{bgk1} for some of these formulas when the
coefficients are complex. It follows from these formulas that any
${\mathfrak R}$-valued function (that is, when $p=q=1$) which is
obtained by addition, multiplication and, when defined,
inversion, of functions of the form \eqref{eq:pol}, is of the
form \eqref{eq:real2}. The matrix-valued case is obtained by
concatenation using the formulas

\[
\begin{split}
\begin{pmatrix}\mathscr H_1&\mathscr H_2\end{pmatrix}(\z)=\\
&\hspace{-2cm}=
\begin{pmatrix}\mathbf D_1&\mathbf D_2\end{pmatrix}+\z
\begin{pmatrix}\mathbf C_1&\mathbf C_2\end{pmatrix}\left(I_{N_1+N_2}-\z
\begin{pmatrix}\mathbf A_1&0\\0&\mathbf
A_2\end{pmatrix}\right)^{-1}
\begin{pmatrix}\mathbf B_1&0\\ 0&\mathbf B_2\end{pmatrix},
\end{split}
\]
and
\[
\begin{split}
\begin{pmatrix}\mathscr H_1\\ \mathscr H_2\end{pmatrix}(\z)=\\
&\hspace{-2cm}=
\begin{pmatrix}\mathbf D_1\\ \mathbf D_2\end{pmatrix}+\z
\begin{pmatrix}\mathbf C_1&0\\ 0&\mathbf C_2
\end{pmatrix}\left(I_{N_1+N_2}-\z
\begin{pmatrix}\mathbf A_1&0\\0&\mathbf
A_2\end{pmatrix}\right)^{-1}
\begin{pmatrix}\mathbf B_1\\ \mathbf B_2\end{pmatrix}
\end{split}
\]
for two functions $\mathscr H_1$ and $\mathscr H_2$ of
appropriate dimensions which admit a realization. Thus, $(1)$
implies $(2)$. We now prove that $(2)$ implies $(3)$. Let
$\mathscr H$ be of the form \eqref{eq:real2}. Then,
\[
R_0^n\mathscr H(\z)=\mathbf C(I_N-\zeta\mathbf A)^{-1}{\mathbf
  A}^{n-1}{\mathbf B},\quad n=1,2,\ldots
\]
But the Cayley-Hamilton theorem holds in any commutative ring
(see for instance \cite[p. 14]{MR839186}, \cite[Theorem 4.3, p.
120]{MR1322960}, \cite[p. A III.107]{MR43:2}). Therefore there
exist an $M\in{\mathbb N}$, a monic polynomial $p$ of degree $M$ with
coefficients in $\mathfrak R$, such that $p({\mathbf A})=0$. It
follows that for $n\ge M$, the function $R_0^n\mathcal H$ is a
linear combination of $1, R_0\mathcal H,\ldots
R_0^{M-1}\mathcal H$ with coefficients in $\mathfrak R$.\\

We now assume that $(3)$ is in force and prove that $(1)$ holds.
First, assume that ${\mathcal H}$ is ${\mathfrak R}$-valued (as opposed
to ${\mathfrak R}^{p\times q}$-valued). By
hypothesis there exists a matrix ${\mathbf A}\in{\mathfrak R}^{M\times
  M}$ such that
\[
R_0\begin{pmatrix} 1&R_0{\mathcal H}&\cdots&R_0^{M-1}\mathcal
H\end{pmatrix} =\begin{pmatrix}1& R_0{\mathcal
H}&\cdots&R_0^{M-1}\mathcal H\end{pmatrix}{\mathbf A}.
\]
Hence
\[
\begin{pmatrix} 1&R_0{\mathcal H}&\cdots&R_0^{M-1}\mathcal
H\end{pmatrix}=\begin{pmatrix}1& R_0{\mathcal
H}&\cdots&R_0^{M-1}\mathcal H\end{pmatrix}(0)(I_M-\zeta{\mathbf
A})^{-1}.
\]
Thus
\[
R_0{\mathcal H}=\mathbf C(I_M-\z{\mathbf A})^{-1}{\mathbf B},
\]
with
\[
\mathbf C=\begin{pmatrix} 1&R_0{\mathcal
H}&\cdots&R_0^{M-1}\mathcal
H\end{pmatrix}(0)\quad{\rm and}\quad {\mathbf B}=\begin{pmatrix}0 \\1\\
  \vdots\\ 0
\end{pmatrix},
\]
and the result follows. The matrix-valued case is treated in much the
same way.\mbox{}\qed\mbox{}\\
\section{Observable pairs}
\setcounter{equation}{0}
Consider an ${\mathfrak R}$-valued function of the complex
variable $\zeta$, of the form \eqref{eq:real2}:
\[
\mathscr H(\z)= \mathbf D+\zeta \mathbf C(I_N-\zeta \mathbf
A)^{-1}\mathbf B.
\]
By Theorem \ref{tm:facto}, we know that $\mathscr H$ is well
defined, in particular for $\zeta$ such that
\[
\det(I_N-\zeta\mathbf A(0))\not=0.
\]
Setting $z=0$ in \eqref{eq:real2} we get the {\sl unperturbed}
transfer function, which motivates Theorem \ref{tm:min} below. We
first give a definition and a proposition.

\begin{Dn}
The pair $({\mathbf C},\mathbf A)\in{\mathfrak R }^{p\times
N}\times{\mathfrak R}^{N\times N}$ is called {\it observable} if the map
\[
\mathbf f\mapsto\begin{pmatrix}\mathbf C\mathbf f&\mathbf C\mathbf
A\mathbf f&\mathbf C\mathbf A^2\mathbf f&\cdots\end{pmatrix}
\]
is injective from ${\mathfrak R}^n$ into $\left({\mathfrak
R}^p\right)^{\mathbb N}$.
\label{dn:obs}
\end{Dn}
See \cite[\S 2.2 p. 58]{MR839186}.\\

Equivalently we have:

\begin{Pn}
\label{pn:obs}
Realization \eqref{eq:real2} is observable if
and only if (with $\mathbf f\in{\mathfrak R}^N$)
\[
\mathbf C(I_N-\z \mathbf A)^{-1}\mathbf f\equiv0_{{\mathfrak
R}}^{p\times N}\Longrightarrow \mathbf f=0_{{\mathfrak R}}^{N}.
\]
\end{Pn}

{\bf Proof:} By Theorem \ref{tm:facto} with $x(t)=(1-\zeta
t)^{-1}$, we have:
\[
(I_N-\zeta\mathbf A)^{-1}=\sum_{p=0}^\infty\z^p{\mathbf A}^p.
\]
Therefore, and using Theorem \ref{300408},
\[
\begin{split}
{\mathbf C}(I_N-\zeta\mathbf A)^{-1}&={\mathbf
C}\sum_{n=0}^\infty\z^p{\mathbf A}^p.
\\
&=\sum_{n=0}^\infty\z^p{\mathbf C}{\mathbf A}^p.
\end{split}
\]
\mbox{}\qed\mbox{}\\

\begin{Tm}
\label{tm:min}
Assume that the realization
\begin{equation}
\label{eq:min0} \mathscr H(\z,0)=\mathbf D(0)+\zeta \mathbf
C(0)(I_N-\zeta \mathbf A(0))^{-1}\mathbf B(0)
\end{equation}
is observable. Then realization \eqref{eq:real2} is observable.
\end{Tm}

{\bf Proof:} Assume first that realization \eqref{eq:min0} is
observable, and let
\[
\mathbf f(z)=\sum_{\alpha\in\ell} f_\alpha z^\alpha, \quad
f_\alpha\in{\mathbb C}^N,
\]
be such that
\begin{equation}
\label{eq:dfgh}
\mathbf C(I_N-\z \mathbf A)^{-1}\mathbf f
%=\sum_{\alpha\in\ell} R_\alpha(\z)H_\alpha(z)
\equiv 0.
\end{equation}
To prove that the realization \eqref{eq:real2} is
observable we need to show that all coefficients $f_\alpha$
in the expansion $\mathbf f(z)=\sum_{\alpha\in\ell} f_\alpha
z^\alpha$, are identically zero.
%where the coefficients $R_\alpha(\z)$ are, as we will see,
%rational functions of $\z$.
Since \eqref{eq:min0} is assumed observable,  setting $z=0$ in
\eqref{eq:dfgh} leads to $f_0=0$. Let us now put
$z=(z_1,z_2,\ldots, z_p,0,0,\ldots)$ in \eqref{eq:dfgh} with $p\ge
1$ and differentiate with respect to $z_1$. We obtain (with
$\mbox{}^\prime$ denoting differentiation with respect to $z_1$)
\begin{equation}
\label{eq:diff}
\begin{split}
\mathbf C^\prime(z)(I_N-\z \mathbf A(z))^{-1}\mathbf f(z)+
\mathbf C(z)\left((I_N-\z \mathbf
A(z))^{-1}\right)^\prime\mathbf f(z)+\\
+ \mathbf C(z)(I_N-\z \mathbf A(z))^{-1}\mathbf
f^\prime(z)\equiv0.
\end{split}
\end{equation}
Setting $z_1=z_2=\cdots=z_p=0$ we obtain that
\[
{\mathbf C}(0)(I_N-\zeta \mathbf
A(0))^{-1}f_{(1,0,0,\ldots)}\equiv0,
\]
and hence $f_{(1,0,0,\ldots)}=0$ since the pair $(\mathbf
C(0),{\mathbf A}(0))$ is observable. Differentiating in
turn \eqref{eq:diff} with respect to $z_1$, we obtain an
expression of the form
\begin{equation}
\label{invalide}
\mathbf X(z)+\mathbf C(z)(I_N-\z\mathbf
A(z))^{-1}f^{\prime\prime}(z)\equiv 0,
\end{equation}
where $\mathbf X$ is a finite sum of the form
\[
\mathbf X(z)=\sum_{j=1}^M
\mathbf U_j(z)f^{(r_j)}(z)
\]
where $r_j\in\left\{0,1\right\}$ and $\mathbf U_j(z)$ is analytic in
$z_1$ and may depend on $\z$. Setting $z_1=0$ in \eqref{invalide}
and taking into account that $f_0=f_{(1,0,0,\ldots)}=0$, we
obtain that
\[
{\mathbf C}(0)(I_N-\zeta \mathbf
A(0))^{-1}f_{(2,0,0,\ldots)}\equiv0,
\]
and hence \( f_{(2,0,0,\cdots)}=0\). More generally, an easy
induction argument shows that the $\alpha_1$-th derivative of
\eqref{eq:dfgh} is of the form
\begin{equation}
\label{invalide1}
\mathbf X(z)+\mathbf C(z)(I_N-\z\mathbf
A(z))^{-1}f^{(\alpha_1)}(z)\equiv 0,
\end{equation}
where $\mathbf X$ is of the form
\[
\mathbf X(z)=\sum_{j=1}^M
\mathbf U_j(z)\mathbf f^{(n_j)}(z),
\]
$\mathbf U_j$ being analytic in $z_1$ and $n_j\in\left\{0,\ldots,
\alpha_1-
1\right\}$.
Setting $z_1=0$ in \eqref{invalide1} we obtain that $f_{(\alpha_1,
0,0,\ldots)}=0$.\\

%Iteratively, we obtain that
%\[
%f_{(\alpha_1,0,0,\ldots)}=0,\quad\forall \alpha_1\in{\mathbb N}.
%\]
Similarly, by setting $z=(0,z_2,0,\ldots)$, and more generally
\[
z=(0,0,0,\ldots, z_j,0,\ldots)
\]
in \eqref{eq:dfgh}, and differentiating, we obtain that
$f_{\alpha}=0$ for all $\alpha\in\ell$ which have only one
non-zero component. We now prove that $f_{(1,1,0,0,\ldots)}=0$.
To that end, set $z=(z_1,z_2,z_3,\cdots, z_p,0,\ldots)$, with
$p\ge 2$, in \eqref{eq:dfgh} and differentiate this equation with
respect to $z_1$ and $z_2$. We obtain an equation of the form
\begin{equation}
\label{sdf}
\mathbf X(z)+{\mathbf C}(z)(I_N-\z \mathbf
A(z))^{-1}\frac{\partial^2\mathbf f}{\partial z_1\partial
z_2}(z)\equiv 0,
\end{equation}
where now $\mathbf X$ is a finite sum of elements of the form
$\mathbf U(z)f(z)$ and $\mathbf U(z)\frac{\partial f}{\partial
z_j} (z)$, with $j\in\left\{1,2\right\}$ and $\mathbf U$ analytic
in $z_1,z_2,\ldots, z_p$. The fact that
\[
f_0=f_{(1,0,0,\ldots)}=f_{(0,1,0,\ldots)}=0\] implies that
$\mathbf X(0)\equiv0$. Setting $z_1=z_2=0$
in \eqref{sdf} then leads to
\[
{\mathbf C}(0)(I_N-\zeta \mathbf
A(0))^{-1}f_{(1,1,0,\ldots)}\equiv0,
\]
and hence $f_{(1,1,0,\ldots)}=0$, where we have used the observability
of the pair $(\mathbf C(0),{\mathbf A}(0))$. By successive
differentiation and setting $z=0$ we obtain that
$f_{(\alpha_1,\alpha_2,0,0,\ldots)}=0$ for every choice of natural
integers $\alpha_1$ and $\alpha_2$. The fact that all other
coordinates $f_\alpha$ are zero, and hence that the pair
$(\mathbf C,\mathbf A)$ is observable, is shown
by induction as follows:\\

{\bf Induction hypothesis:} {\sl For $N\in{\mathbb N}$, it holds that
\begin{equation}
f_{(\alpha_1,\alpha_2,\ldots,\alpha_N,0,0,\ldots)}=0,\quad
\forall (\alpha_1,\alpha_2,\ldots,\alpha_N)\in({\mathbb N}_0)^N,
\end{equation}
and
\begin{equation}
\label{bastille}
\begin{split}
\dfrac{\partial^{\alpha_1+\cdots+\alpha_N}} {\partial
z_1^{\alpha_1}\partial z_2^{\alpha_2}\cdots \partial
z_N^{\alpha_N}} \mathbf
C(z)(I_N-\z \mathbf A(z))^{-1}&=\\
&\hspace{-5cm}= \mathbf X_{\alpha_1,\ldots ,\alpha_N}(z)+\mathbf
C(z)(I_N-\z \mathbf
A(z))^{-1}\dfrac{\partial^{\alpha_1+\cdots+\alpha_N}
  }{\partial z_1^{\alpha_1}\partial z_2^{\alpha_2}\cdots
\partial z_N^{\alpha_N}}\mathbf f(z),
\end{split}
\end{equation}
where
\begin{equation}
z=(z_1,z_2,\ldots , z_p,0,0,\ldots)\quad {\sl with}\quad p\ge N,
\label{eq:z}
\end{equation}
and $\mathbf X_{\alpha_1,\ldots ,\alpha_N}(z)$ is of the form
\begin{equation}
\label{eq:uf}
\mathbf X_{\alpha_1,\ldots ,\alpha_N}(z)=\sum_{j=1}^M
\mathbf U_j(z)\dfrac{\partial^{\beta^{(j)}_1+\cdots+\beta^{(j)}_N}
  }{\partial z_1^{\beta_1^{(j)}}\partial z_2^{\beta_2^{(j)}}\cdots
\partial z_N^{\beta_N^{(j)}}}\mathbf f(z),
\end{equation}
where $\beta^{(j)}_i\le \alpha_i$ for $i=1,\ldots, N$ and
\begin{equation}
\label{ineq}
\beta^{(j)}_1+\cdots+\beta^{(j)}_N<\alpha_1+\cdots+\alpha_N,
\end{equation}
with the functions $\mathbf U_j$ analytic in the variables
$z_1,\ldots , z_p$.}\\

The induction hypothesis holds for $N=1$, as we have shown above.
Assume that it holds at rank $N$.
%We set
%\[
%z=(z_1,z_2,\ldots, z_N,z_{N+1},z_{N+2},\ldots, z_p,0,0,\ldots)
%\]
%in \eqref{eq:dfgh}, with $p\ge N+1$.
We take $p\ge N+1$ in \eqref{eq:z} and differentiate
\eqref{bastille} with respect to $z_{N+1}$.
Since
\[
\begin{split}
\dfrac{\partial}{\partial z_{N+1}}
\mathbf U_j(z)\dfrac{\partial^{\beta^{(j)}_1+\cdots+\beta^{(j)}_N}
  }{\partial z_1^{\beta_1^{(j)}}\partial z_2^{\beta_2^{(j)}}\cdots
\partial z_N^{\beta_N^{(j)}}}\mathbf f(z)&=\\
&\hspace{-20mm}=\left(\dfrac{\partial}{\partial z_{N+1}} \mathbf
U_j(z)\right)\dfrac{\partial^{\beta^{(j)}_1+\cdots+\beta^{(j)}_N}
  }{\partial z_1^{\beta_1^{(j)}}\partial z_2^{\beta_2^{(j)}}\cdots
\partial z_N^{\beta_N^{(j)}}}\mathbf f(z)+\\
\vspace{1mm} \\ &\hspace{-15mm}+ \mathbf
U_j(z)\left(\dfrac{\partial^{\beta^{(j)}_1+\cdots+\beta^{(j)}_N+1}
  }{\partial z_1^{\beta_1^{(j)}}\partial z_2^{\beta_2^{(j)}}\cdots
\partial z_N^{\beta_N^{(j)}}z_{N+1}}\mathbf f(z)\right),
\end{split}
\]
the term
\[
\frac{\partial\mathbf X_{\alpha_1,\ldots ,\alpha_N}(z)}{\partial z_{N+1}}
\]
is of the form
\begin{equation}
\begin{split}
\label{voltaire}
\frac{\partial\mathbf X_{\alpha_1,\ldots ,\alpha_N}(z)}{\partial
  z_{N+1}}
&=\sum_{j=1}^P
\mathbf
V_j(z)\dfrac{\partial^{\beta^{(j)}_1+\cdots+\beta^{(j)}_N+1}
  }{\partial z_1^{\beta_1^{(j)}}\partial z_2^{\beta_2^{(j)}}\cdots
\partial z_N^{\beta_N^{(j)}}\partial z_{N+1}}\mathbf f(z),
\end{split}
\end{equation}
where $P$ is possibly different from $M$ above,
the $\mathbf V_j$ are analytic in the variables $z_1,\ldots,
z_p$, and  the $\beta^{(j)}_i$ are as above. Differentiating the term
\[
\mathbf C(z)(I_N-\z \mathbf
A(z))^{-1}\dfrac{\partial^{\alpha_1+\cdots+\alpha_N}
  }{\partial z_1^{\alpha_1}\partial z_2^{\alpha_2}\cdots
\partial z_N^{\alpha_N}}\mathbf f(z)
\]
in \eqref{bastille} with respect to $z_{N+1}$  we obtain a sum of
two terms. The first is,
\begin{equation}
\label{concorde}
\left(\dfrac{\partial}{\partial z_{N+1}}\mathbf
C(z)(I_N-\z \mathbf
A(z))^{-1}\right)\dfrac{\partial^{\alpha_1+\cdots+\alpha_N}
  }{\partial z_1^{\alpha_1}z_2^{\alpha_2}\cdots
z_N^{\alpha_N}}\mathbf f(z),
\end{equation}
while the second takes the form
\[
\mathbf C(z)(I_N-\z \mathbf
A(z))^{-1}\dfrac{\partial^{\alpha_1+\cdots+\alpha_N+1}
  }{\partial z_1^{\alpha_1}\partial z_2^{\alpha_2}\cdots
\partial z_N^{\alpha_N}z_{N+1}}\mathbf f(z)
\]
This proves \eqref{bastille} for $(\alpha_1,\ldots ,\alpha_N,1)$,
with $\mathbf X_{\alpha_1,\ldots ,\alpha_N,1}$ being the sum of
\eqref{concorde} and of \eqref{voltaire}, that is
\begin{equation}
\label{republique}
\begin{split}
\dfrac{\partial^{\alpha_1+\cdots+\alpha_N+1}} {\partial
z_1^{\alpha_1}
\partial z_2^{\alpha_2}\cdots
\partial z_N^{\alpha_N}
\partial z_{N+1}} \mathbf
C(z)(I_N-\z \mathbf A(z))^{-1}&=\\
&\hspace{-6cm}= \mathbf X_{\alpha_1,\ldots ,\alpha_N,1}(z)+\mathbf
C(z)(I_N-\z \mathbf
A(z))^{-1}\dfrac{\partial^{\alpha_1+\cdots+\alpha_N+1}
  }{\partial z_1^{\alpha_1}\partial z_2^{\alpha_2}\cdots
\partial z_N^{\alpha_N}\partial z_{N+1}}\mathbf f(z).
\end{split}
\end{equation}
Setting $z_1=\cdots =z_p=0$ in that expression, we obtain
\[
\mathbf C(0)(I_N-\z \mathbf A(0))^{-1}f_{(\alpha_1,\ldots,
\alpha_N,1,0,0\ldots)}\equiv0,
\]
and hence $f_{(\alpha_1,\ldots, \alpha_N,1,0,0\ldots)}=0$.
Differentiating \eqref{bastille} a finite number of times  with
respect to $z_{N+1}$, a similar argument will show that
\[
f_{(\alpha_1,\ldots, \alpha_N,\alpha_{N+1},0,0\ldots)}=0,\quad
\forall\alpha_{N+1} \in{\mathbb N}.
\]

%The proof of the corresponding statement for observability is
%proved in the same way, and the claim on minimality then follows.
% We now consider the converse
%statement. We assume that the pair $(\mathbf C,\mathbf A)$ is
%observable, and wish to prove that the pair $(\mathbf
%C(0),\mathbf A(0))$ is also observable. Take in \eqref{eq:dfgh}
%$\mathbf f$ to be a constant.
\mbox{}\qed\mbox{}\\

\begin{Rk}
The converse to the previous theorem does not hold. That is, the
observability of the pair $({\mathbf C}, {\mathbf A})$ does not
imply the observability of the pair $({\mathbf C}(0), {\mathbf
A}(0))$. As an example, take $N=1$ and
\[
\mathbf C(z)=z_1,\quad \mathbf A(z)=1.
\]
The pair $(\mathbf C, \mathbf A)$ is observable since, for
$\mathbf f\in{\mathfrak R}$,
\[
\frac{z_1}{1-\zeta}\mathbf f(z)\equiv0\,\, \Longrightarrow\,\,
\mathbf f=0.
\]
But the pair
\[(\mathbf C(0),\mathbf A(0))=(0,1)\]
is not observable.
\end{Rk}

\section{Controllable pairs and minimal realizations}
\label{5}
\setcounter{equation}{0}
In this section we study controllable pairs and minimal
realizations within the setting of the ring $\mathfrak R$. We
first recall that,  given a commutative ring $R$, one of the
characterization for a pair $( A, B)\in{ R}^{N\times N}\times
{R}^{N\times q}$ of matrices to be controllable (or reachable) is
that the columns of the matrix
\[
\begin{pmatrix}
\ B&A B&\cdots& A^{N-1} B\end{pmatrix}
\]
generate ${R}^N$; see \cite[p. 55]{MR839186}.\\

In the classical case (that is, for the complex numbers, or more
generally, for the case of a field), it is well known that the
pair $(C,A)$ is observable if and only if the pair $(A^T,C^T)$ is
controllable (with $\mbox{}^T$ denoting transpose). This duality
principle does not hold in general in the case of an arbitrary
commutative ring. Only the following direction holds:

\begin{Tm} \cite[Theorem 2.7, p. 59]{MR839186} Let $R$ be a
commutative ring and let $(C,A)\in R^{p\times N}\times R^{N\times
N}$. Assume the pair $(A^T,C^T)$ controllable. Then the pair
$(C,A)$ is observable. \label{tm:bastille}
\end{Tm}

As explained in \cite[p. 59]{MR839186}, the lack of duality comes
form the fact that an homomorphism of modules (say $f$, from the
$R$-module $M_1$ into the $R$-module $M_2$) can be injective
without being residually injective. Recall that residual
injectivity means that, for every maximal ideal $I$ of $R$, the
induced map from $M_1/IR$ into $M_2/IR$ is injective when $f$ is
injective.\\

Theorem \ref{tm:bastille} does not help us to study
controllability based on observability. Furthermore, in \cite[Theorem
2.3 p. 178]{MR1071708}, it is shown that a necessary and
sufficient condition on a commutative ring for the duality
principle to hold for all pairs is that every finitely generated
faithful ideal of the ring contains a unit. As a corollary, the
authors of \cite{MR1071708} state:
\begin{Pn}
\label{pn:opera}
\cite[Corollary 2.4 p. 179]{MR1071708}). If the
duality principle holds in a commutative ring, then the ring is a
total quotient ring.
\end{Pn}
For the purpose of the present paper we do not need to recall the
definition of a faithful ideal (see \cite[Theorem 1.5 (ii), p.
177]{MR1071708}). The total quotient ring of a commutative ring $R$
is the set of formal fractions associated with the set of elements
of $R$ which are not divisors of zero; see
\cite[p. 35]{MR839186}. Thus, in the case of a ring without divisors
of zero
(as is the case for the ring $\mathfrak R$) the total quotient ring is equal
to the quotient field associated with the ring; see \cite[p.
35]{MR839186}. Since ${\mathfrak R}$ is not a field, it follows
that the duality principle is not satisfied on it.\\

After these general preliminaries, let us study controllability
and minimality in the setting of the ring $\mathfrak R$.
 Let us repeat the definition of controllability: The pair
 $({\mathbf A,\mathbf
B})\in{\mathfrak R}^{N\times N}\times {\mathfrak R}^{N\times q}$
is said to be controllable (or reachable) if the columns of the
matrix
\[
\begin{pmatrix}
\mathbf B&\mathbf A\mathbf B&\cdots&\mathbf A^{N-1}\mathbf
B\end{pmatrix}
\]
generate ${\mathfrak R}^N$. See \cite[p. 55]{MR839186}. We
therefore have:

\begin{Pn}
\label{opera}
Assume the pair $({\mathbf A,\mathbf
B})\in{\mathfrak R}^{N\times N}\times {\mathfrak R}^{N\times q}$
to be controllable. Then the pair $({\mathbf A(0),\mathbf
B}(0))\in{\mathbb C}^{N\times N}\times {\mathbb C }^{N\times q}$
is controllable.
\end{Pn}

{\bf Proof:} Since ${\mathbb C}^N\subset {\mathfrak R}^N$, for
every $f\in {\mathbb C}^N$ there exists $\mathbf a\in{\mathbf
R}^{Nq}$ such that
\[
f=\begin{pmatrix}\mathbf B&\mathbf A\mathbf B&\cdots&\mathbf
A^{N-1}\mathbf B\end{pmatrix}(z)\mathbf a(z).
\]
Setting $z=0$ in this equality we get the controllability of the
pair $({\mathbf A(0),\mathbf B}(0))$.\mbox{}\qed\mbox{}\\

\begin{Rk} We note the difference between Theorem \ref{tm:min} and
Proposition \ref{opera}. In the former, observability at $z=0$
implies observability in ${\mathfrak R}$. In the latter,
controllability in ${\mathfrak R}$ implies controllability at
$z=0$.
\end{Rk}

The converse of Proposition \ref{opera} would be an analogue of
Theorem \ref{tm:min} for the case of controllable pairs. But this
is not possible for the ring ${\mathfrak R}$, in view of
Proposition \ref{pn:opera}, since $\mathfrak R$ is different from
its total quotient ring (which is in fact its quotient field since
$\mathfrak R$ has no divisors of zero).\\

Still, we can give a counterpart of Theorem \ref{tm:min} for
controllable and minimal realizations with the following ad-hoc
definitions:

\begin{Dn}
\label{dn:min1} Realization \eqref{eq:real2} will be called
${\mathfrak R}$-controllable if the following condition holds: Let
$\mathbf f\in{\mathfrak R}^{1\times N}$. Then:
\[
\mathbf f(I_N-\z \mathbf A)^{-1}\mathbf B\equiv0_{{\mathfrak
R}}^{1\times q}\Longrightarrow \mathbf f=0_{{\mathfrak
R}}^{1\times N}.
\]
A realization will be called ${\mathfrak R}$-minimal if it is both
observable and $\mathfrak R$-controllable.
\end{Dn}

We can then state:

\begin{Tm}
\label{tm:min2} Assume that realization \eqref{eq:min0}
\begin{equation*}
\mathscr H(\z,0)=\mathbf D(0)+\zeta \mathbf C(0)(I_N-\zeta
\mathbf A(0))^{-1}\mathbf B(0)
\end{equation*}
is controllable (resp. minimal). Then realization
\eqref{eq:real2} is  $\mathfrak R$-controllable (resp.
 $\mathfrak R$-minimal).
\end{Tm}

{\bf Proof:} The first statement is proved as
Theorem \ref{tm:min}. The second statement follows then from the
definition of minimality.
\mbox{}\qed\mbox{}\\

\section{Hilbert-space valued transfer functions}
\setcounter{equation}{0} We recall that the Hilbert spaces
${\mathcal H}_{k}$ have been defined above by the finiteness of
the norm \eqref{michelle}. By $\mathbf I(\mathcal H_k)$ we mean
the image of $\mathcal H_k$ under the Hermite transform. We note
that
\[
{\mathfrak R}=\cup_{k=1}^\infty \mathbf I(\mathcal H_k).
\]
\begin{Tm}
Let $\mathcal H$ be given by realization \eqref{eq:real2},
and let $l,k$ be natural numbers such that $k>l+1$. Assume that,
in the state space equations \eqref{eq:sp2}, the entries of $ A$
and $ C$ are in ${\mathcal H}_{l}$ and the entries of $ B$ and ${
D}$ are in ${\mathcal
  H}_{k}$.
Then, the transfer function ${\mathscr H}$ is $\mathbf I({\mathcal
H}_{k})$-valued.
\end{Tm}

{\bf Proof:}
Inequality \eqref{vage} expresses the fact that
the multiplication operator
\[
T_h\,:\,u\mapsto h\lozenge u
\]
is a bounded map from the Hilbert space ${\mathcal H}_{k}$ into
itself. Therefore the entries of the ${\mathfrak R}^{p\times
q}$-valued function $ \mathbf C{\mathbf A}^n{\mathbf B}$ are in
$\mathbf I({\mathcal H}_{k})$. To conclude the proof, it remains to show
that for every complex number $\z$ such that $(I_N-\z {\mathbf
A})$ is invertible, the power series
\[
\sum_{n=0}^\infty \z^n\mathbf C{\mathbf A}^n{\mathbf B}\]
converges in $\mathbf I({\mathcal H}_{k})$ to ${\mathscr H}(\z)$. But this
is a consequence of Theorem \ref{tm:facto}.
\mbox{}\qed\mbox{}\\

Using once more V\r{a}ge's inequality \eqref{vage} we have:
\begin{Cy}
Let now $m>k+1$, where $k$ is as in the previous theorem.
Then, the operator of multiplication by ${\mathscr H}$ sends
$\mathbf I({\mathcal H}_{m})$-valued signals into
$\mathbf I({\mathcal H}_{m})$-valued signals.
\end{Cy}

\bibliographystyle{plain}
\def\cprime{$'$} \def\lfhook#1{\setbox0=\hbox{#1}{\ooalign{\hidewidth
  \lower1.5ex\hbox{'}\hidewidth\crcr\unhbox0}}} \def\cprime{$'$}
  \def\cprime{$'$} \def\cprime{$'$} \def\cprime{$'$} \def\cprime{$'$}

%\bibliography{/users/faculty/math/dany/bib/all}
%\bibliography{all}
\end{document}